\documentclass[12pt]{article}

\usepackage{geometry} 
\usepackage{latexsym}
\usepackage{amssymb}
\usepackage{amsthm}
\usepackage{amsmath,amsfonts}
\usepackage{xcolor}
\usepackage{float}

\def\titlerunning#1{\gdef\titrun{#1}}
\makeatletter
\def\author#1{\gdef\autrun{\def\and{\unskip, }#1}\gdef\@author{#1}}
\def\address#1{{\def\and{\\\hspace*{18pt}}\renewcommand{\thefootnote}{}%
\footnote {#1}}%
\markboth{\autrun}{\titrun}}
\makeatother
\def\email#1{\hspace*{4pt}{\em e-mail}: #1}

\newtheorem{thm}{Theorem}[section]

\theoremstyle{definition}
\newtheorem{rem}[thm]{Remark}

\titlerunning{}

\title{The existence of some directed strongly regular graphs on 54 and 108 vertices}
\author{Dean Crnkovi\' c, Andrea \v Svob and Matea Zubovi\'c \v Zutolija}

\begin{document}
\maketitle

\address{D. Crnkovi\' c, A. \v Svob, M. Zubovi\'c \v Zutolija: Faculty of Mathematics, University of Rijeka, Radmile Matej\v ci\'c 2, 51000 Rijeka, Croatia;
\email{\{deanc,asvob,matea.zubovic\}@math.uniri.hr}
}

\begin{abstract}
In this paper, we prove the existence of directed strongly regular graphs with parameters (108,11,3,2,1), (108,14,10,0,2), (108,22,12,6,4), (108,23,9,8,4), (108,25,15,8,5), (108,34,18,12,10), (108,38,22,12,14), (108,39,23,14,14), (108,41,35,16,15), (108,42,33,18,15) and (108,46,22,19,20). Further, we obtain directed strongly regular graphs with parameters (54,10,4,1,2), (54,11,4,3,2) and (54,16,7,4,5). The constructions are obtained by considering finite groups acting transitively on 54 and 108 vertices. 
\end{abstract}

\bigskip



\section{Introduction}

A directed graph is an ordered pair $\Delta = (\mathcal{V},\mathcal{E})$ consisting of a vertex set $\mathcal{V}$ and a set $\mathcal{E}$ of ordered pairs of vertices called arcs. If $(x,y) \in \mathcal{E}$, we write $x \to y$.

The concept of a directed strongly regular graph was introduced by Duval \cite{Duval88} as a natural directed version of strongly regular graphs. Specifically, a directed strongly regular graph (dsrg) with parameters $(v,k,t,\lambda,\mu)$ is a directed graph on $v$ vertices in which every vertex has both indegree and outdegree equal to $k$, and in which the number of directed paths of length two from $x$ to $y$ equals $t$ if $x=y$, equals $\lambda$ if $x \to y$ is an arc, and equals $\mu$ otherwise. Throughout this paper we restrict attention to digraphs without loops and without multiple arcs between the same ordered pair of vertices.

Given two directed graphs $\Delta_1=(\mathcal{V}_1,\mathcal{E}_1)$ and $\Delta_2=(\mathcal{V}_2,\mathcal{E}_2)$, a bijection $f:\mathcal{V}_1 \to \mathcal{V}_2$ is called an \textit{isomorphism} if
\[
(x,y) \in \mathcal{E}_1 \iff (f(x),f(y)) \in \mathcal{E}_2
\]
for all $x,y \in \mathcal{V}_1$. An isomorphism from $\Delta$ to itself is an \textit{automorphism}. The collection of all automorphisms of $\Delta$ forms the full automorphism group $Aut(\Delta)$, and any of its subgroups is referred to as an automorphism group of $\Delta$.

Unlike strongly regular graphs, whose existence for small parameter sets is largely settled, the existence question for directed strongly regular graphs remains open for a substantial number of feasible parameter sets, even on relatively small vertex sets.

Recently, in \cite{bcs, BCSZ} some new directed strongly regular graphs were constructed using simple groups. In this paper, we applied the method introduced in \cite[Theorem 3]{cms} to the finite groups, not necessarily simple groups. We demonstrate that the method is effective also when applied beyond the simple groups setting.

The feasibility table maintained in \cite{BH} lists a number of parameter sets on $108$ vertices for which no directed strongly regular graph had yet been found. The present note closes these gaps by giving first examples of dsrgs with parameters $(108,11,3,2,1)$, $(108,14,10,0,2)$, $(108,22,12,6,4)$, $(108,23,9,8,4)$, $(108,25,15,8,5)$, $(108,34,18,12,10)$, $(108,38,22,12,14)$, $(108,39,23,14,14)$, $(108,41,35,16,15)$, $(108,42,33,18,15)$ and $(108,46,22,19,20)$. For each of the eleven parameter sets above we constructed directed strongly regular graphs, by studying transitive actions of finite groups on $108$ points. In total, we constructed 45 mutually nonisomorphic directed strongly regular graphs having 108 vertices. Moreover, by studying transitive actions of finite groups on $54$ points we obtained directed strongly regular graphs with parameters (54,10,4,1,2) and (54,16,7,4,5). As a consequence of finding new directed strongly regular graphs with parameters (108,23,9,8,4),  we prove the existence of directed strongly regular graphs with parameters (54,11,4,3,2), which are related to the directed strongly regular graphs with parameters (108,23,9,8,4) by the construction called blowup (\cite{Duval88}).

We note that these examples are of interest from a group-theoretic standpoint: the underlying graphs are relatively small, yet several of the groups producing them act with large rank.

\bigskip

The computations in this paper are made by using programs written for Magma \cite{magma}. The directed strongly regular graphs constructed in this paper can be found at the link:

\begin{verbatim}
https://www.math.uniri.hr/~matea.zubovic/DSRGs_108vertices/dsrgs_108.html
\end{verbatim}

\section{Construction of new directed strongly regular graphs}

We constructed the directed strongly regular graphs by using the method described in \cite[Theorem 3]{cms}.
That construction produces simple 1-designs on which a group $G$ acts transitively on the points and blocks. Hence, if the incidence structure of a 1-design obtained using \cite[Theorem 3]{cms} is the adjacency matrix of a directed strongly regular graph, then the graph admits a transitive action of $G$ on the set of vertices. In this paper, we construct the directed strongly regular graphs using the groups $AGL(2,3), S_3\wr S_3, 3^3.S_4$ and $6^2:(2\times S_3)$.

\subsection{Construction of dsrgs from the group $AGL(2,3)$}

The group $AGL(2,3)$ (it is the SmallGroup(432,734) in the Small Groups Library) is a finite group of order 432, and has, up to conjugation, exactly 46 subgroups. Among them, there are exactly two subgroups of order 4, one subgroup $H_1$, isomorphic to the cyclic group of order $4$, having the property that $AGL(2,3)$ acts on cosets of $H_1$ in 32 orbits, and another subgroup $H_2$, isomorphic to the group $2^2$, having the property that $AGL(2,3)$ acts on cosets of $H_2$ in 33 orbits.

Using the method from \cite[Theorem 3]{cms}, by taking $G=AGL(2,3)$ and the stabilizer of a vertex $G_{\alpha}=H_1$, we constructed two nonisomorphic dsrgs with parameters (108,14,10,0,2) and two nonisomorphic dsrgs with parameters (108,22,12,6,4) admitting a transitive action of the group $AGL(2,3)$. Precisely, we constructed two nonisomorphic dsrgs(108,14,10,0,2), which we denote by $\Delta_1$ and $\Delta_2$, where one digraph is obtained from the other by reversing the arcs. We say that the digraphs $\Delta_1$ and $\Delta_2$ are reverse to each other (see \cite[Section 3.4]{BCSZ}). The full automorphism group of the digraphs $\Delta_1$ and $\Delta_2$ is isomorphic to $AGL(2,3)$. 

Further, by taking $G=AGL(2,3)$ and the stabilizer of a vertex $G_{\alpha}=H_1$, we constructed two nonisomorphic dsrgs(108,22,12,6,4), which we denote by $\Delta_3$ and $\Delta_4$, where one digraph is obtained from the other by reversing the arcs. The full automorphism group of the digraphs $\Delta_3$ and $\Delta_4$ is isomorphic to $AGL(2,3)$.

By taking $G=AGL(2,3)$ and the stabilizer of a vertex $G_{\alpha}=H_2$, we constructed six different parameter sets for $k=11,22,25,34,41,46$ as follows.

We obtained four nonisomorphic dsrgs with parameters (108,11,3,2,1), which we denote by $\Delta_5$, $\Delta_6$, $\Delta_7$ and $\Delta_8$. The digraph $\Delta_5$ can be obtained by reversing the arcs of $\Delta_6$, and vice versa. The same holds for digraphs $\Delta_7$ and $\Delta_8$. The full automorphism group of the digraphs $\Delta_5$ and $\Delta_6$ is isomorphic to $AGL(2,3)$, and the full automorphism groups of the digraphs $\Delta_7$ and $\Delta_8$ are isomorphic to $(3^4.2)\wr (A_4.2)$.

We obtained two nonisomorphic dsrgs with parameters (108,22,12,6,4), which we denote by $\Delta_9$ and $\Delta_{10}$, which are reverse to each other. The full automorphism group of the digraphs $\Delta_9$ and $\Delta_{10}$ is isomorphic to $AGL(2,3)$. We note here that digraphs $\Delta_9$ and $\Delta_{10}$ are not isomorphic to $\Delta_3$ and $\Delta_4$.

We obtained two nonisomorphic dsrgs with parameters (108,25,15,8,5), which we denote by $\Delta_{11}$ and $\Delta_{12}$, which are reverse to each other. The full automorphism group of the digraphs $\Delta_{11}$ and $\Delta_{12}$ is isomorphic to $AGL(2,3)$.

We obtained two nonisomorphic dsrgs with parameters (108,34,18,12,10), which we denote by $\Delta_{13}$ and $\Delta_{14}$, which are reverse to each other. The full automorphism group of the digraphs $\Delta_{13}$ and $\Delta_{14}$ is isomorphic to $AGL(2,3)$.

We obtained two nonisomorphic dsrgs with parameters (108,41,35,16,15), which we denote by $\Delta_{15}$ and $\Delta_{16}$, which are reverse to each other. The full automorphism group of the digraphs $\Delta_{15}$ and $\Delta_{16}$ is isomorphic to $ASL(2,3).S_3$.

We obtained 16 nonisomorphic dsrgs with parameters (108,46,22,19,20), which we denote by $\Delta_{17}, \dots, \Delta_{32}$. The constructed digraphs have the full automorphism group isomorphic to $AGL(2,3)$. Among these 16 nonisomorphic digraphs, we obtained eight pairs of reversed digraphs.

\subsection{Construction of dsrgs from the group $S_3\wr S_3$}

The group $S_3\wr S_3$ (the SmallGroup(1296,3490)) is a finite group of order 1296, and has, up to conjugation, exactly 21 subgroups of index 108. Among them, there is exactly one subgroup of order 12, $H_3$, that is isomorphic to $D_6$ and having the property that $S_3\wr S_3$ acts on cosets of $H_3$ in 19 orbits. By taking the subgroup $H_3$ as the stabilizer of a vertex, we obtained two nonisomorphic 
dsrgs (108,38,22,12,14) admitting a transitive action of the group $S_3\wr S_3$, which we denote by $\Delta_{33}$ and $\Delta_{34}$. The constructed digraphs have the full automorphism groups isomorphic to $S_3\wr S_3$. Constructed dsrgs are reverse to each other.

\subsection{Construction of dsrgs from the group $3^3.S_4$}

The group $3^3.S_4$ (it is the SmallGroup(648,703)) is a finite group of order 648, and has, up to conjugation, exactly 8 subgroups of index 108. Among them, there is exactly one subgroup of order 6, $H_4$, isomorphic to $S_3$ and having the property that $3^3.S_4$ acts on cosets of $H_4$ in 22 orbits. By taking the subgroup $H_4$ as the stabilizer of a vertex, we obtained two nonisomorphic 
dsrgs (108,42,33,18,15) admitting a transitive action of the group $3^3.S_4$, which we denote by $\Delta_{35}$ and $\Delta_{36}$. The constructed digraphs have the full automorphism groups isomorphic to $3^3.S_4$. Constructed dsrgs are reverse to each other.

\subsection{Construction of dsrgs from the group $6^2:(2\times S_3)$}\label{108-54}

The group $6^2:(2\times S_3)$ (the SmallGroup(432,523)) is a finite group of order 432, and has, up to conjugation, exactly 9 subgroups of order 4. Among them, there are exactly two subgroups, $H_5$ and $H_6$, isomorphic to the group $2^2$, having the property that the group $6^2:(2\times S_3)$ acts on cosets of these subgroups in 40 orbits. Both subgroups yield the same results.

Using the method from \cite[Theorem 3]{cms}, by taking $G=6^2:(2\times S_3)$ and the stabilizer of a vertex $G_{\alpha}=H_5$ (or $G_{\alpha}=H_6$), we constructed five nonisomorphic dsrgs with parameters (108,23,9,8,4) admitting a transitive action of the group $6^2:(2\times S_3)$, which we denote by $\Delta_{37}, \Delta_{38}, \Delta_{39}, \Delta_{40}, \Delta_{41}$. The full automorphism group of the constructed digraphs is isomorphic to $2^{36}.3.2^6.2^6.A_4^2.2^3.2$. Among these nonisomorphic digraphs, we obtained two pairs of reversed digraphs and one self-reversed digraph.

Moreover, the group $6^2:(2\times S_3)$ has, up to conjugation, exactly one subgroup of order 4, $H_7$, isomorphic to the group $2^2$, having the property that the group $6^2:(2\times S_3)$ acts on cosets of these subgroups in 33 orbits. Using the subgroup $H_7$, we constructed four nonisomorphic dsrgs with parameters (108,39,23,14,14) admitting a transitive action of the group $6^2:(2\times S_3)$, which we denote by $\Delta_{42}, \Delta_{43}, \Delta_{44}, \Delta_{45}$. The full automorphism groups of the digraphs $\Delta_{42}$ and $\Delta_{43}$ are isomorphic to the group
$(2^2 \times ((2^4 : (3^2 : 3)) : 2)) : 2$ and of $\Delta_{44}$ and $\Delta_{45}$ are isomorphic to the group $6^2:(2\times S_3)$. Digraphs $\Delta_{42}$ and $\Delta_{43}$, and $\Delta_{44}$ and $\Delta_
{45}$, are reverse to each other.

\subsubsection{Construction of dsrgs with parameters (54,11,4,3,2)}

In \cite{Duval88}, the following method for constructing directed strongly regular graphs was given. We refer to \cite{BCSZ} as additional literature in which the technique was applied also.

Given a directed strongly regular graph $D$ with
parameters $(n,k,t,\lambda,\mu)$, where $t = \mu$, and an integer $m$,
one can replace each vertex by a set of $m$ mutually nonadjacent vertices,
all with the same in-neighbours and out-neighbours,
and in that way obtain a directed strongly regular graph $D'$ with parameters $(mn,mk,mt,m\lambda,m\mu)$. We call $D'$ the $m$-blowup of $D$.

Similarly, if we apply the previous construction to the complementary graph, we get the following. From a directed strongly regular graph $D$ with
parameters $(n,k,t,\allowbreak\lambda,\mu)$, where $t = \lambda+1$,
one obtains a directed strongly regular graph $D'$ with parameters
$(mn,m(k+1)-1,m(t+1)-1,m(\lambda+2)-2,m\mu)$ by replacing
each vertex by a set of $m$ mutually adjacent vertices.

Taking the five nonisomorphic dsrgs with parameters (108,23,9,8,4) admitting a transitive action of the group $6^2:(2\times S_3)$, i.e. $\Delta_{37}, \Delta_{38}, \Delta_{39}, \Delta_{40}, \Delta_{41}$ (constructed in \ref{108-54}), one can check that they can be seen as a 2-blowups of five nonisomorphic dsrgs with parameters (54,11,4,3,2), which we denote by $\Delta_{46},\dots,\Delta_{50}$.

\begin{rem}
The five nonisomorphic dsrgs with parameters (54,11,4,3,2), whose
2-blowups are the five nonisomorphic DSRGs with parameters
$(108,23,9,8,4)$, can be constructed using the method from \cite[Theorem 3]{cms}. In Table \ref{new-54} we give details how one can construct them. Among them, there are two pairs of reversed digraphs and one self-reversed digraph.
\end{rem}

\subsubsection{Construction of dsrgs with parameters (54,10,4,1,2) and (54,16,7,4,5)}

The group $6^2:(2\times S_3)$ (the SmallGroup(432,523)) is a finite group of order 432, and has, up to conjugation, exactly 7 subgroups of order 8. Among them, there is exactly one subgroup, namely $H_8$, isomorphic to the group $2^3$, having the property that the group $6^2:(2\times S_3)$ acts on cosets of these subgroups in 30 orbits.

Using the method from \cite[Theorem 3]{cms}, by taking $G=6^2:(2\times S_3)$ and the stabilizer of a vertex $G_{\alpha}=H_8$, we constructed 18 nonisomorphic dsrgs with parameters (54,16,7,4,5) admitting a transitive action of the group $6^2:(2\times S_3)$, which we denote by $\Delta_{51}, \dots, \Delta_{68}$. The full automorphism group of the constructed digraphs is isomorphic to $3.S_3^2$. Among these nonisomorphic digraphs, we obtained nine pairs of reversed digraphs.

Furthermore, by taking $G=6^2:(2\times S_3)$ and the stabilizer of a vertex $G_{\alpha}=H_8$, we constructed two nonisomorphic dsrgs with parameters (54,10,4,1,2), admitting a transitive action of the group $6^2:(2\times S_3)$, which we denote by $\Delta_{69}, \Delta_{70}$. The full automorphism group of the constructed digraphs is isomorphic to $3.S_3^2$. Two digraphs are reverse to each other.

\section{Overview of Results}

It follows that, in every case where the number of nonisomorphic graphs is even, no graph is isomorphic to its reverse. However, when the number of nonisomorphic graphs is odd, exactly one graph is isomorphic to its reverse. A summary of these results is provided in Tables \ref{new-108} and \ref{new-54}.

\begin{table}[H]
\begin{footnotesize}
\centering
\begin{tabular}{|c|c|c|c|c|c|c|}
\hline
$G$ & $H$ & rank & $(n,k,t,\lambda,\mu)$ & $\#$ nonisom. & Aut$\,\Delta$\\
	\hline
$AGL(2,3)$ &4 & 32& (108,14,10,0,2) & 2 & $AGL(2,3)$  \\
$AGL(2,3)$ &4 & 32& (108,22,12,6,4) & 2 & $AGL(2,3)$  \\
$AGL(2,3)$ &$2^2$ & 33& (108,11,3,2,1) & 4 & $AGL(2,3)$ (2),  \\
 && & &  & $(3^4.2)\wr (A_4.2)$ (2)  \\
$AGL(2,3)$ &$2^2$ & 33& (108,22,12,6,4) & 2 & $AGL(2,3)$   \\
$AGL(2,3)$ &$2^2$ & 33& (108,25,15,8,5) & 2 & $AGL(2,3)$   \\
$AGL(2,3)$ &$2^2$ & 33& (108,34,18,12,10) & 2 & $AGL(2,3)$  \\
$AGL(2,3)$ &$2^2$ & 33& (108,41,35,16,15) & 2 & $ASL(2,3).S_3$  \\
$AGL(2,3)$ &$2^2$ & 33& (108,46,22,19,20) & 16 & $AGL(2,3)$  \\ \hline

$S_3\wr S_3$ &$D_6$ & 19& (108,38,22,12,14) & 2 & $S_3\wr S_3$\\ \hline

$3^3.S_4$ &$S_3$ & 22& (108,42,33,18,15) & 2 & $3^3.S_4$ \\ \hline

$6^2:(2\times S_3)$ &$2^2$ & 40& (108,23,9,8,4) & 5 & $2^{36}.3.2^6.2^6.A_4^2.2^3.2$ \\ 
$6^2:(2\times S_3)$ &$2^2$ & 33& (108,39,23,14,14) & 4 & $(2^2 \times ((2^4:(3^2:3)):2)):2$ (2), \\ 
 & & & & & $6^2:(2\times S_3)$ (2)  \\ \hline
\hline
\end{tabular}
\caption{New DSRGs having 108 vertices}\label{new-108}
\end{footnotesize}
\end{table}

In Table \ref{new-54}, we give directed strongly regular graphs having 54 vertices.

\begin{table}[H]
\begin{footnotesize}
\centering
\begin{tabular}{|c|c|c|c|c|c|c|}
\hline
$G$ & $H$ & rank & $(n,k,t,\lambda,\mu)$ & $\#$ nonisom. & Aut$\,\Delta$\\
	\hline
$6^2:(2\times S_3)$& $2^3$ & 30& (54,11,4,3,2) & 5 & $3.S_3^2$  \\
$6^2:(2\times S_3)$& $2^3$ & 30& (54,10,4,1,2) & 2 &  $3.S_3^2$  \\
$6^2:(2\times S_3)$& $2^3$ & 30& (54,16,7,4,5) & 18 & $3.S_3^2$  \\
\hline
\end{tabular}
\caption{New DSRGs having 54 vertices}\label{new-54}
\end{footnotesize}
\end{table}

\section{Statements and Declarations}

\subsection{Declaration of competing interests}

The authors declare no conflict of interest.

\subsection{Data Availability Statement}
The directed strongly regular graphs constructed in this paper can be found at the link:

\begin{verbatim}
https://www.math.uniri.hr/~matea.zubovic/DSRGs_108vertices/dsrgs_108.html
\end{verbatim}

\section*{Acknowledgement}
The authors would like to thank A. E. Brouwer for pointing out the construction of dsrgs on 54 vertices and for verifying all the results, which significantly improved the quality of the paper. This work was supported by the Croatian Science Foundation under the project number HRZZ-IP-2022-10-4571 and by European Union-NextGenerationEU, project number uniri-iz-25-46-KonGeoGraGru.

\end{document}